\documentclass{amsart}
\usepackage{amssymb,latexsym,enumerate}
\newtheorem{Thm}{Theorem}
\theoremstyle{remark}
\newtheorem*{Ack*}{Acknowledgement}
\DeclareMathOperator{\GL}{GL}
\begin{document}

\title[Right orderable residually finite $p$-groups]
{Right orderable residually finite $p$-groups and a Kourovka notebook
problem}

\author[P. A. Linnell]{Peter A. Linnell}

\address{Math \\ VPI \\ Blacksburg \\ VA 24061--0123
\\ USA}

\email{linnell@math.vt.edu}
\urladdr{http://www.math.vt.edu/people/linnell/}

\begin{abstract}
A.~H.~Rhemtulla proved that if a group is a residually finite
$p$-group for infinitely many primes $p$, then it is two-sided
orderable.  In problem 10.30 of the Kourovka notebook 14\,th.\
edition, N.~Ya.~Medvedev asked if there is a non-right-orderable
group which is a residually finite $p$-group for at least two
different primes $p$.  Using a result of Dave Witte, we will show that
many subgroups of finite index in $\GL_3(\mathbb {Z})$ give examples
of such groups.  On the other hand we will show that no such example
can exist among solvable by finite groups.
\end{abstract}

\keywords{right orderable group, residually finite $p$-group}

\subjclass{Primary: 20F60; Secondary: 06F15}

\maketitle

A group $G$ is right orderable if it has a total ordering $\le $ such
that $x \le y \Rightarrow xg \le yg$ whenever $g,x,y \in G$.  It is
\emph{two-sided orderable} if in addition $x \le y \Rightarrow gx
\le gy$ whenever $g,x,y \in G$.  For much
information on right ordered groups, see the books
\cite{KopytovMedvedev96,MuraRhemtulla77}.  It was proved in
\cite{Rhemtulla73} that a group which is a residually finite
$p$-group for infinitely many primes $p$ is two-sided orderable.  In
problem 10.30 of the Kourovka notebook 14\,th.\ edition
\cite{Kourovka99}, N.~Ya.~Medvedev asks if there is a
non-right-orderable group which
is a residually finite $p$-group for at least two different primes
$p$.  We shall prove
\begin{Thm} \label{Texample}
Let $\mathcal {P}$ be a finite set of primes.  Then there exists a
non-right-orderable group which is a residually finite $p$-group for
all $p \in \mathcal {P}$.
\end{Thm}
The proof depends on a theorem of Witte \cite{Witte94}, and we will
see that the groups in Theorem \ref{Texample} can be taken to be
subgroups of finite index in $\GL_3(\mathbb {Z})$, the group of 3
by 3 invertible matrices with integer entries.

Going in the opposite direction, we shall prove the
following result.
\begin{Thm} \label{Torder}
Let $p,q$ be distinct primes and let $G$ be a solvable by finite
group.  If $G$ is a residually
finite $p$-group and a residually
finite $q$-group, then $G$ has a series
\[
1 = G_0 \lhd \dots \lhd G_n = G
\]
with $G_i \lhd G$ and $G_{i+1}/G_i$ torsion free abelian for all
$i$.  In particular $G$ is right orderable.
\end{Thm}

\begin{Ack*}
This work was carried out while I was at the
Sonderforschungsbereich in M\"unster.  I would like to thank
Wolfgang L\"uck for organizing my visit to M\"unster, and the
Sonderforschungsbereich for financial support.
\end{Ack*}

\begin{proof}[Proof of Theorem \ref{Texample}]
For each prime $p$, let $G_p = \{A \in \GL_3(\mathbb {Z}) \mid A
\equiv I \mod p \}$, the congruence subgroup of level $p$ (here $I$
denotes the identity matrix of $\GL_3(\mathbb {Z})$).  Then $G_p$ is
a residually finite $p$-group.  Set $G = \bigcap_{p \in \mathcal
{P}} G_p$.  Then $G$ is a residually finite $p$-group for all $p \in
\mathcal {P}$.  Furthermore $G$ has finite index in $\GL_3(\mathbb
{Z})$, so by \cite[Proposition 3.3]{Witte94} we see that $G$ is not
right orderable.  This completes the proof of Theorem
\ref{Texample}.
\end{proof}

\begin{proof}[Proof of Theorem \ref{Torder}]
Observe that $G$ is torsion free.
Let $H$ be a normal solvable subgroup of finite index in $G$ with
minimal derived length.  We shall prove the result by induction on
the derived length of $H$, the result being obvious when this is zero
because then $H=1$ and thus $G$ will be finite.  We now assume that
the derived length of $H$ is at least one.

Let $A_0$ denote the penultimate term of the derived series of $H$.
Then $A_0$ is a normal abelian subgroup of $G$ and $H/A_0$ has
strictly smaller derived length than
that of $H$.  Let $A$ be a maximal
normal abelian subgroup of $G$ containing $A_0$.

We shall let $\hat{K}$ denote the pro-$p$ completion of a group
$K$.  Then the exact sequence $1 \to A \to G \to
G/A \to 1$ yields an exact sequence
\[
\hat {A} \longrightarrow \hat{G} \longrightarrow
\widehat{G/A} \longrightarrow 1.
\]
Let $B$ denote the image of $\hat {A}$ in $\hat{G}$.
Since $A$ is abelian, we see that $\hat{A}$ is
abelian and we deduce that $B$ is abelian.  Also we
may view $G$ as a subgroup of $\hat{G}$ because $G$ is a
residually finite $p$-group.  Therefore $B \cap G$ is an abelian
normal subgroup of $G$ containing $A$
and we conclude that $B \cap G = A$.  But
$B \cap G$ can also be described as the kernel
of the natural map $G \to
\widehat{G/A}$.  Therefore $G/A$ is isomorphic to a
subgroup of $\widehat{G/A}$ and we deduce that
$G/A$ is a residually finite $p$-group.  Similarly
$G/A$ is a residually finite $q$-group.  Since the derived
length of $HA/A$ is strictly less than the derived
length of $H$, induction shows that $G$ has a series
\[
A = G_1 \lhd \dots \lhd G_n = G
\]
with $G_i \lhd G$ and $G_{i+1}/G_i$ torsion free abelian for all
$i$.  By setting $G_0 = 1$, we obtain the required series for the
first part of Theorem \ref{Torder}.

The assertion that $G$ is right orderable now follows from
\cite[Theorem 7.3.2]{MuraRhemtulla77} and the fact that torsion free
abelian groups are right orderable.
\end{proof}

\bibliographystyle{plain}

\begin{thebibliography}{1}

\bibitem{KopytovMedvedev96}
Valeri{\u\i}~M. Kopytov and Nikola{\u\i}~Ya. Medvedev.
\newblock {\em Right-ordered groups}.
\newblock Consultants Bureau, New York, 1996.

\bibitem{Kourovka99}
V.~D. Mazurov and E.~I. Khukhro, editors.
\newblock {\em The {K}ourovka notebook}.
\newblock Russian Academy of Sciences Siberian Division Institute of
  Mathematics, Novosibirsk, augmented edition, 1999.
\newblock Unsolved problems in group theory.

\bibitem{MuraRhemtulla77}
Roberta~Botto Mura and Akbar Rhemtulla.
\newblock {\em Orderable groups}.
\newblock Marcel Dekker Inc., New York, 1977.
\newblock Lecture Notes in Pure and Applied Mathematics, Vol. 27.

\bibitem{Rhemtulla73}
A.~H. Rhemtulla.
\newblock Residually ${F}\sb{p}$-groups, for many primes $p$, are orderable.
\newblock {\em Proc. Amer. Math. Soc.}, 41:31--33, 1973.

\bibitem{Witte94}
Dave Witte.
\newblock Arithmetic groups of higher {$\mathbb {Q}$}-rank cannot act on
  $1$-manifolds.
\newblock {\em Proc. Amer. Math. Soc.}, 122(2):333--340, 1994.

\end{thebibliography}

\end{document}